\ifx\documentclass\undefined
   \documentstyle[12pt,a4wide,twoside]{article}
\else
   \documentclass[12pt,twoside]{article}
   \usepackage{a4wide}
\fi

\usepackage{pstricks}
\usepackage{psfrag}
\usepackage{amssymb}

\usepackage{graphicx}

\usepackage{color}

\newcommand {\nc}   {\newcommand}
\nc {\be}   {\begin{equation}}
\nc {\ee}   {\end{equation}}
\nc {\beq}  {\begin{eqnarray}}
\nc {\eeq}  {\end{eqnarray}}
\nc {\beqs} {\begin{eqnarray*}}
\nc {\eeqs} {\end{eqnarray*}}

\nc {\supp} {{\rm supp}}
\nc {\diag} {{\rm diag}}

\nc{\D}     {\displaystyle}
\nc{\SSS}   {\scriptscriptstyle}

\nc{\comment} [1] {}

\newcount\hour \newcount\minute
\hour=\time    \divide\hour by 60
\minute=\hour  \multiply\minute by 60 \advance\minute by -\time \minute=-\minute
\def\twodigits#1{\ifnum #1 < 10{0#1}\else{#1}\fi}

\renewcommand{\comment}[1]{}

\font\Blackbrd=msbm10 scaled 1200       

\newcommand{\jump}[1]   {\mbox{$\big[ \hspace{-0.7mm} \big[ #1
            \big] \hspace{-0.7mm} \big]$} }

\newcommand{\sumt}  {\sum_{T\in\caT}}

\newcommand{\osc}{\mathop{\rm osc}\nolimits}

\newcommand{\caE}{{\mathcal E}}

\newcommand{\caN}{{\mathcal N}}

\newcommand{\caT}{{\mathcal T}}

\renewcommand{\le}  {\lesssim}      

\newcommand{\R}     {\mbox{\Blackbrd R}}    
\newcommand{\Poly}  {\mbox{\Blackbrd P}}    

\renewcommand{\div} {{\rm{div} \,}}     

\newtheorem{theorem}{Theorem}[section]
\newtheorem{lemma}[theorem]{Lemma}
\newtheorem{corollary}[theorem]{Corollary}
\newcommand{\bt}{\begin{theorem}}
\newcommand{\et}{\end{theorem}}
\newcommand{\br}{\begin{remark}}
\newcommand{\er}{\end{remark}}
\newcommand{\bc}{\begin{corollary}}
\newcommand{\ec}{\end{corollary}}
\newcommand{\bl}{\begin{lemma}}
\newcommand{\el}{\end{lemma}}
\newcommand{\bp}{\begin{proposition}}
\newcommand{\ep}{\end{proposition}}
\newcommand{\bd}{\begin{definition}}
\newcommand{\ed}{\end{definition}}
\newcommand{\bex}{\begin{example}}
\newcommand{\eex}{\end{example}}

\newtheorem{definition}[theorem]{Definition}
\newtheorem{example}[theorem]{Example}

\newtheorem{remark}[theorem]{Remark}

\newenvironment{proof}{
{\noindent \bf Proof:}}{\quad \hfill \rule{2mm}{2mm}\medskip}

\author{Serge Nicaise and Juliette Venel \footnote{Universit\'e de
Valenciennes et du Hainaut Cambr\'esis, LAMAV, FR CNRS 2956,
Institut des Sciences et Techniques de Valenciennes, F-59313 -
Valenciennes Cedex 9 France, email:
Serge.Nicaise,Juliette.Venel@univ-valenciennes.fr}}

\begin{document}

  \title{A posteriori error estimates  for a
   finite element approximation of  transmission problems with sign changing coefficients}

  {
\maketitle
 \noindent

{
    \begin{abstract}
 We perform the a posteriori error analysis of   residual type of a
transmission problem with sign changing coefficients. According to
\cite{BonnetCiarlet09} if the  contrast is large enough, the
continuous problem can be transformed into a coercive one.   We
further show that a similar property holds for the discrete problem
for any regular
  meshes, extending the framework from \cite{BonnetCiarlet09}. The reliability and
efficiency of the proposed estimator is  confirmed by some numerical
tests.
\end{abstract}



 }

    \thispagestyle{plain}

\section{Introduction}
\label{introduction} Recent years have witnessed a growing interest
in the study of diffusion problems with a sign changing coefficient.
These problems appear in several areas of physics, for example in
electromagnetism
\cite{Engheta:02,MaWolf:95,Maystre:04,Pendry:00,Ramdani99}. Thus
some mathematical investigations have been performed and concern
existence results \cite{BonnetDauge99,Ramdani99} and numerical
approximations by the finite element methods
\cite{Ramdani99,BonnetCiarlet06,BonnetCiarlet08,BonnetCiarlet09},
with some a priori error analyses. But for such problems the
regularity of the solution may be poor and/or unknown and
consequently an a posteriori error analysis would be more
appropriate. This analysis is the aim of the present paper.

For continuous Galerkin finite element methods, there now exists a
large amount of literature on a posteriori error estimations for
(positive definite) problems in mechanics or electromagnetism.
Usually   locally defined a posteriori error estimators are
designed. We refer the reader to the monographs
\cite{AinsworthOden,BS01,Monk03,verfurth:96b} for a good overview on
this topic.

In contrast to the recent paper \cite{BonnetCiarlet09} we will not
use quasi-uniform meshes that are not realistic for an a posteriori
error analysis.
That is why we improve their finite element analysis in order to
allow only regular meshes in Ciarlet's sense \cite{ciarlet:78}.

The paper is structured as follows:  We recall in Section 2 the
"diffusion" problem and the technique from
\cite{BonnetCiarlet09} that allows to establish its well-posedness
for sufficiently large contrast.
In Section 3, we prove that the discrete approximation is well-posed
by introducing an ad-hoc discrete lifting operator. The a posteriori error analysis is
performed in Section 4, where upper and lower bounds are obtained.
Finally in Section 5 some numerical tests are presented that confirm
the reliability and efficiency of our estimator.

Let us finish this introduction with some notations used in the
remainder of the paper: On $D$, the $L^2(D)$-norm will be denoted by
$\|\cdot\|_D$.  The usual norm and semi-norm of $H^{s}(D)$ ($s\geq
0$) are denoted by $\|\cdot\|_{s,D}$ and $|\cdot|_{s,D}$,
respectively. In  the case $D=\Omega$, the index $\Omega$ will be
omitted.
Finally, the notations $a\lesssim b$ and $a\sim b$ mean the
existence of positive constants $C_1$ and $C_2$, which are
independent of the mesh size and of the considered quantities $a$
and $b$ such that $a\leq C_2b$ and $C_1b\leq a\leq C_2b$,
respectively. In other words, the constants may depend on the aspect
ratio of the mesh and the diffusion coefficient (see below).

\section{The boundary value problem}
\label{sbvp}

Let $\Omega$ be a bounded open domain of $\R^2$  with boundary
$\Gamma$. We suppose that $\Omega$  is split up into two sub-domains
$\Omega_+$ and $\Omega_-$ with a Lipschitz boundary  that we suppose
to be polygonal
in such a way
that
\[\bar \Omega=\bar \Omega_+\cup \bar\Omega_-, \quad \Omega_+\cap
\Omega_-=\emptyset,
\]
see Figure \ref{fig1} for an example.

\begin{figure}
\centering
\begin{pspicture}(2,2)
\psline(1,2)(0,2)(0,0)(2,0)(2,1) \psline(1,1)(1,2)(2,2)(2,1)(1,1)
\rput(1.5,1.5){$\Omega_+$} \rput(0.7,0.5){$\Omega_-$}
\end{pspicture}
 \caption{The domain $\Omega$
\label{fig1}}
\end{figure}
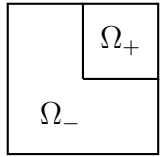

\noindent We now assume that the diffusion coefficient $a$ belongs
to $L^\infty(\Omega)$ and is positive (resp. negative) on $\Omega_+$
(resp. $\Omega_-$). Namely there exists $\epsilon_0>0$ such that
\beq\label{ass_a+} a(x)\geq \epsilon_0, \textmd{ for} \mbox{ a. e. }
x\in
\Omega_+,\\
\label{ass_a-} a(x)\leq -\epsilon_0, \textmd{ for} \mbox{ a. e. } x\in
\Omega_-.\eeq

\noindent In this situation we consider the following   second order boundary
value problem with Dirichlet boundary conditions: \be
\label{divgradu} \left\{\begin{array}{rcll}
-\div(a\hbox{~}\nabla u) &=&f &\mbox{~in~$\Omega$},\\
u&=&0 &\hbox{~on~$\Gamma$}.
\end{array}\right.
\ee

\noindent The variational formulation of (\ref{divgradu}) involves the
bilinear form
$$ {\it B}(u,v)=\int_\Omega   a \nabla u \cdot \nabla v    $$
and the Hilbert space
\[
H^1_0(\Omega)=\{u\in H^1(\Omega): u=0 \hbox{~on~}\Gamma\}.
\]

Due to the lack of coercivity of ${\it B}$ on $H^1_0(\Omega)$ (see
\cite{BonnetDauge99,BonnetCiarlet08,BonnetCiarlet09}), this problem
does not fit into a standard framework. In
\cite{BonnetCiarlet08,BonnetCiarlet09}, the proposed approach is to
use a bijective and continuous linear mapping $\mathbb T$ from
$H^1_0(\Omega)$ into itself that allows to come back to the coercive
framework. Namely these authors assume that ${\it B}(u, {\mathbb T}
v)$ is coercive in the sense that there exists $\alpha>0$ such that
\be\label{tcoercivite} {\it B}(u, {\mathbb T} u)\geq \alpha
\|u\|_{1,\Omega}^2 \quad \forall u\in H^1_0(\Omega). \ee Hence given
$f \in L^2(\Omega)$, by the Lax-Milgram theorem the problem \be
\label{FV} {\it B}(u,{\mathbb T} v)=\int_\Omega f {\mathbb T} v
\quad \forall v \in H^1_0(\Omega), \ee has a unique solution $u \in
H^1_0(\Omega)$. Since ${\mathbb T}$  is an isomorphism,  the
original problem \be \label{FV} {\it B}(u,v)=\int_\Omega f v \quad
\forall v \in H^1_0(\Omega), \ee has also a unique solution $u \in
H^1_0(\Omega)$.

In \cite{BonnetCiarlet09}, the mapping ${\mathbb T}$ is built by
using a trace lifting operator $\cal R$ from $H^{1/2}_{00}(\Sigma)$
into $H^1_-(\Omega_-)$, where $\Sigma=\partial \Omega_-\cap \partial
\Omega_+$ is the interface between $\Omega_-$ and $\Omega_+$,
\[
H^1_\pm(\Omega_\pm)=\{u\in H^1(\Omega_\pm): u=0 \mbox{ on } \partial
\Omega_\pm\setminus \Sigma\},
\]
and
\[
H^{1/2}_{00}(\Sigma)=\{u_{|\Sigma}: u\in H^1_-(\Omega_-)\}=
\{u_{|\Sigma}: u\in H^1_+(\Omega_+)\}
 \] is the
space of the restrictions to $\Sigma$ of functions in
$H^1_-(\Omega_-)$ (or in $H^1_+(\Omega_+)$). This last space may be
equipped with the norms
\[ \displaystyle
\|p\|_{1/2, \pm}=\inf_{u\in H^1_\pm(\Omega_\pm)\atop p=u_{|\Sigma}}
 |u|_{1,\Omega_\pm}.
\]
With the help of such a lifting, a possible mapping ${\mathbb T}$ is
given by (see \cite{BonnetCiarlet09})
\[
{\mathbb T} v= \left\{\begin{array}{ll}
v_+   \mbox{~in~}  \Omega_+,\\
-v_-+2{\cal R} (v_{+|\Sigma}) \hbox{~in~} \Omega_-,
\end{array}\right.
\]
where $v_\pm$ denotes the restriction of $v$ to $\Omega_\pm$. With
this choice, it is shown in Proposition 3.1 of
\cite{BonnetCiarlet09} that (\ref{tcoercivite}) holds if
\be\label{star} K_{\cal R}=\sup_{v\in H^1_+(\Omega_+)\atop v\ne
0}\frac{|{\it B}_-({\cal R} (v_{|\Sigma}),{\cal R}
(v_{|\Sigma}))|}{{\it B}_+(v,v)} <1, \ee where ${\it B}_\pm(u,v)=
\displaystyle \int_{\Omega_\pm}   a \nabla u \cdot \nabla
v$. \\

For concrete applications, one can make the following particular
choice for ${\cal R}$, that we denote by ${\cal R}_p$: for any
$\varphi\in H^{1/2}_{00}(\Sigma)$ we define ${\cal R}_p(\varphi)=w$
as the unique solution $w\in H^1_-(\Omega_-)$ of
\[
\Delta w= 0 \hbox{ in } \Omega_-, \quad w=\varphi \hbox { on }
\Sigma.
\]
With this choice, one obtains that $K_{{\cal R}_p}<1$ if the contrast
$$
\frac{\min_{\Omega_-} |a|}{\max_{\Omega_+} a}
$$ is large enough, we refer to Section 3 of \cite{BonnetCiarlet09}
for more details.

\br Note that in \cite{BonnetCiarlet08,BonnetCiarlet09} the authors
consider sub-domains $\Omega_+$ and $ \Omega_-$ with a
pseudo-Lipschitz boundary. However the previous arguments from
 \cite{BonnetCiarlet09} (shortly summarized above)
 are not valid in this case since the space $H^1(\Omega_+)$ equipped with the  norm
 $|\cdot|_{1,\Omega_+}$ is not complete. \er

\section{The discrete  approximated problem}
\label{s_approx}

Here we consider the following standard Galerkin approximation of
our continuous
 problem.
We consider a triangulation $\caT$ of $\Omega$, that is a
"partition" of $\Omega$ made of triangles $T$ (closed subsets of
$\bar \Omega$) whose edges are denoted by $e$. We
 assume that this triangulation is regular, i.e.,
 for any element $T$, the ratio $h_{T}/\rho _{T}$ is bounded by a constant $\sigma >0$
 independent of $T$ and of the   mesh size $h=\max_{T\in \caT} h_T$,
 where $h_{T}$ is the diameter of $T$ and $\rho _{T}$ the diameter of its largest inscribed ball.
 We further assume that $\caT$ is conforming with the partition   of $\Omega$, i.e.,
 each triangle is assumed to be either included into   $\bar \Omega_+$
 or into   $\bar \Omega_-$.
 With each edge  $e$ of the triangulation,
  we denote by $h_e$ its length and $n_{e}$  a unit normal vector
  (whose orientation can be arbitrary chosen) and the so-called patch  $\omega_e=\cup_{e\subset T}
  T$,
  the union of triangles having $e$ as edge. We similarly associate with each vertex $x$, a patch $\omega_x =\cup_{x\in T}
  T$.
For a triangle $T$, $n_T$
  stands for the outer unit normal vector of $T$.
  $\caE$ (resp. $\caN$) represents the set of edges (resp. vertices) of the triangulation.
In the sequel, we need to distinguish between edges (or vertices)  included into
$\Omega$ or into $\Gamma$, in other words, we set
\beqs \caE_{int} &=&\{e\in\caE: e\subset \Omega\},\\
\caE_{\Gamma} &=&\{e\in\caE: e\subset  \Gamma\},\\
\caN_{int} &=&\{x\in\caN: x\in \Omega\}.
\eeqs

\noindent Problem (\ref{FV}) is approximated by  the  continuous  finite
element  space: \begin{equation}\label{defxh} V_h=\left\{v_h \in
H^1_0(\Omega): v_{h|T} \in \Poly_\ell(T),\, \forall T \in \caT\right\},
\end{equation}
where $\ell$ is a fixed positive integer and the space $\Poly_\ell(T) $ consists of polynomials of degree at most $\ell$.

\noindent The Galerkin approximation of problem
 (\ref{FV}) reads now: Find $u_h \in V_h$,
  such that
\be\label{FVdiscrete} {\it B}(u _{h},v _{h}) =\int_\Omega f v_h
\quad \forall v_{h} \in V_h.\ee

\noindent Since there is no reason that the bilinear form would be coercive on
$V_h$, as in \cite{BonnetCiarlet09} we need to use a discrete
mapping ${\mathbb T}_h$ from $V_h$ into itself defined by  (see
\cite{BonnetCiarlet09})
\[
{\mathbb T}_h v_h= \left\{\begin{array}{ll}
v_{h+}   \mbox{~in~}  \Omega_+,\\
-v_{h-}+2{\cal R}_h (v_{h+|\Sigma}) \hbox{~in~} \Omega_-,
\end{array}\right.
\]
where ${\cal R}_h$ is a discrete version of the operator ${\cal R}$.
Here contrary to \cite{BonnetCiarlet09} and in order to avoid the
use of quasi-uniform meshes (meaningless in an a posteriori error
analysis), we take \be\label{defRh} {\cal R}_h=I_h{\cal R}, \ee
where $I_h$ is a sort of Cl\'ement interpolation operator
\cite{clement:75} and ${\cal R}$ is any  trace lifting operator
 from $H^{1/2}_{00}(\Sigma)$ into $H^1_-(\Omega_-)$ (see the previous section). More
precisely for $\varphi_h\in H_h(\Sigma)=\{v_{h|\Sigma}: v_h\in
V_h\}$, we set
\[
 I_h{\cal R}(\varphi_h)=\sum_{x\in \caN_-} \alpha_x \lambda_x,
\]
where $\caN_-=\caN_{int}\cap \bar \Omega_-$, $\lambda_x$ is the
standard hat function (defined by $\lambda_x\in V_h$ and satisfying
$\lambda_x(y)=\delta_{xy}$) and $\alpha_x\in \R$ are defined by
\[
\alpha_x= \left\{\begin{array}{ll}
|\omega_x|^{-1}\int_{\omega_x} {\cal R}(\varphi_h) &  \mbox{~if~}  x\in \caN_{int}\cap \Omega_-, \vspace{6pt}\\
\varphi_h(x) & \hbox{~if~} x\in \caN_{int}\cap \Sigma,
\end{array}\right.
\]
where we recall that $\omega_x$ is the patch associated with $x$,
which is simply the support of $\lambda_x$. Note that $I_h$
coincides with the Cl\'ement interpolation operator $I_{\rm Cl}$ for
the nodes in $\Omega_-$ and only differs on the nodes on $\Sigma$.
Indeed let us recall the definition of  $I_{\rm Cl}{\cal
R}(\varphi_h)$ (defined in a Scott-Zhang manner
\cite{LRScott_SZhang_1992a} for the points belonging to $\Sigma$):
$$  I_{\rm Cl}{\cal R}(\varphi_h)=\sum_{x\in \caN_-} \beta_x \lambda_x$$
with
\[
\beta_x= \left\{\begin{array}{ll}
|\omega_x|^{-1}\int_{\omega_x} {\cal R}(\varphi_h) &  \mbox{~if~}  x\in \caN_{int}\cap \Omega_-, \vspace{6pt}\\
|e_{x}|^{-1}\int_{e_{x}}  {\cal R}(\varphi_h) d\sigma & \hbox{~if~}
x\in \caN_{int}\cap \Sigma \hbox{ with } e_{x}=\omega_x\cap \Sigma.
\end{array}\right.
\]
The definition of $ I_h$ aims at ensuring that
\[
I_h{\cal R}(\varphi_h)= \varphi_h \hbox{ on } \Sigma.
\]

\noindent Let us now prove that ${\cal R}_h$ is uniformly bounded.
\bt\label{tRhborne} For all $h>0$ and $\varphi_h \in H_h(\Sigma)$,
one has
\[
|{\cal R}_h(\varphi_h)|_{1,\Omega_-}\lesssim \| \varphi_h\|_{1/2,-}.
\]
\et
\begin{proof}
For the sake of simplicity we make the proof in the case $\ell=1$,
the general case is treated in the same manner by using modified
Cl\'ement interpolation operator.

\noindent  Since ${\cal R}$ is bounded from
 $H^{1/2}_{00}(\Sigma)$ into   $H^1_-(\Omega_-)$, one has
\be\label{borneRp} |{\cal R}(\varphi_h)|_{1,\Omega_-}\lesssim \|
\varphi_h\|_{1/2,-}. \ee Hence it suffices to show that
\be\label{Rhborne} |(I-I_h){\cal R}(\varphi_h)|_{1,\Omega_-}\lesssim
\| \varphi_h\|_{1/2,-}. \ee
 For that purpose, we distinguish the triangles $T$ that have no nodes
 in $\caN_{int}\cap \Sigma$ to the other ones:
 \\

 1. If $T$ has no nodes in $\caN_{int}\cap \Sigma$, then $I_h
 {\cal R}(\varphi_h)$ coincides with $I_{\rm Cl} {\cal
 R}(\varphi_h)$ on $T$ and therefore by a standard property of the
 Cl\'ement interpolation operator, we have
\be \label{Rhborne1}
 |(I-I_h){\cal R}(\varphi_h)|_{1,T}=|(I-I_{\rm Cl}){\cal
 R}(\varphi_h)|_{1,T}\lesssim \| {\cal
 R}(\varphi_h)\|_{1,\omega_T},
\ee where the patch $\omega_T$ is
 given by $\omega_T=\displaystyle \bigcup_{T'\cap T\ne\emptyset} T'$.
 \\

 2. If $T$ has at least one node in $\caN_{int}\cap \Sigma$, by the
 triangle inequality we may write
 \[
 |(I-I_h){\cal R}(\varphi_h)|_{1,T}\leq |(I-I_{\rm Cl}){\cal
 R}(\varphi_h)|_{1,T}+ |(I_{\rm Cl}-I_h){\cal
 R}(\varphi_h)|_{1,T}.
 \]
  For the first term of this right-hand
 side we can still use (\ref{Rhborne1}) and therefore it remains to
 estimate the second term.
For that one, we notice that
$$(I_{\rm Cl}-I_h){\cal R}(\varphi_h)=\sum_{x\in T\cap \Sigma}
 (\alpha_x-\beta_x)\lambda_x  \hbox{ on } T.
$$
Hence $$  |(I_{\rm Cl}-I_h){\cal R}(\varphi_h)  |_{1,T} \lesssim
\sum_{x\in T\cap \Sigma}  |\alpha_x-\beta_x|.$$ Since ${\cal
R}(\varphi_h)=\varphi_h$ on $\Sigma$ and due to the definition of
$I_{\rm Cl}$, it follows that for $x \in T \cap \Sigma$,
$$  |\alpha_x-\beta_x  |= \left | \varphi_h(x)- |e_{x}|^{-1}\int_{e_{x}}
  \varphi_h d\sigma\right |. $$
Since all norms are equivalent in finite dimensional spaces, we have
for all $v_h \in \Poly_1(e_x)$, \be | v_h(x)| \lesssim |e_x|^{-1/2}
\|v_h \|_{e_x} .  \label{equivnorm}\ee Moreover, \be |e_x|^{-1/2}
\left\| \varphi_h- |e_x|^{-1}  \int_{e_{x}} \varphi_h  d\sigma
\right\|_{e_x} \lesssim   | \varphi_h|_{1/2,e_x},
\label{seminorm}\ee
 where here $|\cdot|_{1/2, e_x}$ means the standard
 $H^{1/2}(e_x)$-seminorm.
Thus Inequalities (\ref{equivnorm}) with $ v_h =\varphi_h-
|e_{x}|^{-1}\int_{e_{x}}
  \varphi_h d\sigma$ and (\ref{seminorm}) imply that
$$ | \alpha_x-\beta_x | \lesssim   | \varphi_h|_{1/2,e_x} .$$

\noindent  All together we have shown that \be\label{Rhborne2}
|(I-I_h){\cal
 R}(\varphi_h)|_{1,T}  \lesssim \| {\cal
 R}(\varphi_h)\|_{1,\omega_T\cap \bar \Omega_-}+|\varphi_h |_{1/2,
 \omega_T\cap \Sigma}. \ee

 \noindent Taking the sum of the square of (\ref{Rhborne1}) and of
 (\ref{Rhborne2}), we obtain that \[ |(I-I_h){\cal R}(\varphi_h)|_{1,\Omega_-}^2 \lesssim \| {\cal
 R}(\varphi_h)\|_{1, \Omega_-}^2+|\varphi_h |_{1/2, \Sigma}^2. \]
 We conclude thanks to (\ref{borneRp}) and to the fact that
 \[
 |\varphi_h |_{1/2, \Sigma}\lesssim \|\varphi_h \|_{1/2,-}.
\]
\end{proof}

 \noindent This Theorem and Proposition 4.2 of \cite{BonnetCiarlet09}
 allow to conclude that (\ref{FVdiscrete}) has a unique solution
 provided that
(\ref{star}) holds, in particular if
 the contrast is large
 enough.

 \noindent Note that the advantage of our construction of ${\cal R}_h$ is that
we no more need the quasi-uniform property of the meshes imposed in
\cite{BonnetCiarlet09}.

\section{The  a posteriori error analysis}

Error estimators can be constructed in many different ways as, for
example, using residual type error estimators which measure locally
the jump of the discrete flux \cite{verfurth:96b}. A different
method, based on equilibrated fluxes, consists in solving local
Neumann boundary value problems \cite{AinsworthOden} or in using
Raviart-Thomas interpolant
\cite{ainsworth:06,cochez:07,ern_nic:07,ern:07}. Here since the
coercivity constant is not explicitly known, we chose the simplest
approach of residual type.

The residual estimators
are denoted by \beq \label{defestiR}
\eta_{R}^2=\sum_{T \in \caT} \eta_{R,T}^2, \quad \eta_{J}^2=\sum_{T \in
\caT} \eta_{J,T}^2, \eeq where the indicators $\eta_{R,T}$ and
$\eta_{J,T}$ are defined by \beqs \eta_{R,T}&=&h_T\|f_T+\div(a\nabla u_h)\|_T,\\
\eta_{J,T}&=&\sum_{e\in \caE_{int}: e\subset
T}h_e^{1/2}\|\jump{a\nabla u_h\cdot n_e} \|_e,\eeqs when $f_T$ is an
approximation of $f$, for instance
\[
 f_T=|T|^{-1}\int_T f.
 \]
Note that $\eta_{R,T}^2$ is meaningful if $a_{|T}\in W^{1,1}(T),$
for all $T\in \caT$.

\subsection{Upper bound} \label{s_upperRT1}

\bt \label{tupperbound} Assume that $a\in L^\infty(\Omega)$
satisfies (\ref{ass_a+})-(\ref{ass_a-})
 and that $a_{|T}\in W^{1,1}(T),$ for all $T\in
\caT$. Assume further that (\ref{star}) holds. Let $u\in
H^1_0(\Omega)$ be the unique solution of Problem (\ref{FV}) and let
$u_h$ be its Galerkin approximation, i.e. $u_h \in V_h$ a solution
of (\ref{FVdiscrete}). Then  one has \be \label{upperbound}
\|\nabla(u-u_h)\|\lesssim \eta_{R}+\eta_J+\osc(f), \ee
 where
 \[
 \osc(f)=\left(\sumt h_T^2 \|f-f_T\|^2\right)^\frac12.
\]
 \et
\begin{proof}
By the coerciveness assumption (\ref{tcoercivite}), we may write
  \be\label{decerror}
\|\nabla ( u-u_h)\|^2\lesssim {\it B}(u-u_h, {\mathbb T} (u-u_h)).
\ee But we notice that the Galerkin relation \[ {\it B}(u-u_h,v_h)=0
\quad \forall v_h\in V_h\] holds.
 Hence by taking $v_h=I_{Cl} {\mathbb T} (u-u_h)$, (\ref{decerror})
 may be written
\be\label{decerror2} \|\nabla ( u-u_h)\|^2\lesssim {\it B}(u-u_h,
(I-I_{Cl}){\mathbb T} (u-u_h)). \ee Now we apply standard arguments,
see for instance \cite{verfurth:96b}. Namely applying element-wise
Green's formula and writing for shortness $w=(I-I_{Cl}){\mathbb T}
(u-u_h)$, we get \beqs \|\nabla ( u-u_h)\|^2\lesssim -\sumt
\int_T\div(a\nabla (u-u_h)) w \\
+\sum_{e\in \caE_{int}} \int_e \jump{a\nabla (u-u_h)\cdot n} w\, d\sigma,
\eeqs reminding that $w=0$ on $\Gamma$.
  By Cauchy-Schwarz's inequality we directly obtain
\beqs \|\nabla ( u-u_h)\|^2\lesssim \sumt
\|f+\div(a\nabla u_h)\|_T \|w\|_T \\
+\sum_{e\in \caE_{int}} \|\jump{a\nabla u_h\cdot n}\|_e \|w\|_e.
\eeqs By standard interpolation error estimates, we get \beqs
\|\nabla (u-u_h)\|^2\lesssim \Big(\sumt h_T^2 \|f+\div(a\nabla
u_h)\|_T^2
\\
+\sum_{e\in \caE_{int}} h_e\|\jump{a\nabla u_h\cdot
n}\|_e^2\Big)^{1/2} |{\mathbb T} (u-u_h)|_{1,\Omega}. \eeqs Since
${\mathbb T}$ is an isomorphism, we conclude that \beqs\|\nabla (
u-u_h)\| \lesssim \Big(\sumt h_T^2 \|f+\div(a\nabla u_h)\|_T^2
\\
+\sum_{e\in \caE_{int}} h_e\|\jump{a\nabla u_h\cdot
n}\|_e^2\Big)^{1/2}. \eeqs This
leads to the conclusion due to the triangle inequality.
\end{proof}

\subsection{Lower bound} \label{s_lowerRT1} The  lower bound is
fully standard since by a careful reading of the proof of
Proposition 1.5 of \cite{verfurth:96b}, we see that it does not use
the positiveness of the diffusion coefficient $a$. Hence we can
state the

  \bt\label{tlowerbound} Let the assumptions of Theorems
\ref{tupperbound}   be satisfied. Assume furthermore that $a_{|T}$
is constant for all $T\in \caT$. Then for each element $T \in \caT$
the following estimate holds
\[ \eta_{R,T}+\eta_{J,T}  \le
|u-u_h |_{1,\omega_T}+\osc(f,\omega_T), \] where
\beqs\osc(f,\omega_T)^2=\sum_{T'\subset\omega_T}h_{T'}^2\|f-f_T'\|_{T'}^2.
\eeqs \et

\section{Numerical results}
\label{secnum}
\subsection{The polynomial solution}
\label{poly} In order to illustrate our theoretical predictions,
this first numerical test consists in validating our computations on
a simple case, using an uniform refinement process. Let $\Omega$ be
the square $(-1,1)^2$, $\Omega_+=(0,1)\times (-1,1)$ and
$\Omega_-=(-1,0)\times (-1,1)$. We assume that ${a}=1$ on $\Omega_+$
and $a=\mu<0$ on $\Omega_-$. In such a situation we can take
\[
{\cal R}(v_+)(x,y)=v_+(-x,y) \quad \forall (x,y)\in \Omega_-.
\]
With this choice we see that
\[
K_{\cal R}=|\mu|,
\]
and therefore for $|\mu|<1$, (\ref{tcoercivite}) holds and Problem
(\ref{FV}) has a unique solution. We further easily check that the
corresponding mapping $\mathbb T$ is an isomorphism since  $(\mathbb
T)^2=\mathbb T$. Similarly by exchanging the role of $\Omega_+$ and
$\Omega_-$, (\ref{tcoercivite}) will also hold if $|\mu|>1$.

\noindent Now we take as exact solution $$\begin{array}{ll}
u(x,y)=\mu x(x+1)(x-1)(y+1)(y-1) & \forall (x,y)\in \Omega_+,\\
u(x,y)=x (x+1)(x-1)(y+1)(y-1) & \forall (x,y)\in
\Omega_-,  \end{array}$$
$f$ being fixed accordingly.

Let us recall that $u_h$ is the finite element solution, and set
$e_{L^2}(u_h)=\| u-u_h\|$ and $e_{H^1}(u_h)=\|u-u_h\|_1$ the $L^2 $
and $H^1 $ errors. Moreover let us define $\eta(u_h)=
\eta_{R}+\eta_J$ the estimator and  $CV_{L^2}$  (resp. $CV_{H^1}$)
as the experimental convergence rate of the error $e_{L^2}(u_h)$
(resp. $e_{H^1}(u_h)$) with respect to the mesh size defined by
$DoF^{-1/2}$, where the number of degrees of freedom is $DoF$,
computed from one line of the table to the following one.

\noindent Computations are performed with $\mu=-3$ using a  global
mesh refinement process from an initial cartesian grid. First, it
can be seen from Table \ref{tablepoly} that the convergence rate of
the $H^1$ error norm is equal to one, as theoretically expected (see
\cite{BonnetCiarlet09}). Furthermore the convergence rate of the
$L^2 $ error norm is 2, which is a consequence of the Aubin-Nitsche
trick and regularity results for Problem  (\ref{divgradu}). Finally,
the reliability of the estimator is ensured since the ratio in the
last column (the so-called  effectivity index), converges towards a
constant close to 6.5.
\begin{table}[ht]
\centering
\begin{tabular}{|l|l|c|c|c|c|c|}
\hline
&&&&&& \\
$k$ & $DoF$ & $e_{L^2}(u_h)$ & $CV_{L^2}$ & $ e_{H^1}(u_h)$ & $CV_{H^1}$ & $\displaystyle \frac{\eta(u_h)  }{e_{H^1}(u_h)  }$  \\
&&& &&&\\
\hline
1 & 289  & 2.37E-02 & & 5.33E-01& &6.70 \\
\hline
2 & 1089 & 5.95E-03 &2.08& 2.67E-01&1.04& 6.59 \\
\hline
3 & 4225 & 1.49E-03 &2.04& 1.34E-01& 1.02 & 6.53  \\
\hline
4 & 16641 & 3.73E-04 &2.02& 6.68E-02&1.01& 6.49 \\
\hline
5 & 32761 &1.89E-04 &2.01& 4.75E-02&1.01&6.48\\
\hline
6 & 90601 &6.79E-05 &2.01& 2.85E-02& 1.00 & 6.47\\
\hline
7 &  251001& 2.45E-05  & 2.00 & 1.71E-02 &1.00 & 6.47  \\
\hline
\end{tabular}
\caption{The polynomial solution with $\mu=-3$  (uniform refinement).}
\label{tablepoly}
\end{table}

\subsection{A singular solution}

Here we analyze an example introduced in \cite{BonnetDauge99} and
precise some results from \cite{BonnetDauge99}. The domain
$\Omega=(-1,1)^2$ is decomposed into two sub-domains $\Omega_+=(0,1)
\times (0,1)$, and $\Omega_-=\Omega\setminus \bar\Omega_+$, see
Figure \ref{fig1}. As before we take
 $a=1$ on $\Omega_+$ and $a=\mu<0$
 on $\Omega_-$.
According to Section 3 of \cite{BonnetDauge99}, Problem (\ref{FV})
has a singularity $S$ at $(0,0)$ if $\mu<-3$ or if $\mu\in (-1/3,0)$
given in polar coordinates by \beqs S_+(r,\theta)=r^\lambda
(c_1\sin(\lambda \theta)+ c_2 \sin(\lambda (\frac{\pi}{2}-\theta)))
&&\hbox{ for } 0<\theta<\frac{\pi}{2},\\
S_-(r,\theta)=r^\lambda (d_1\sin(\lambda (\theta-\frac{\pi}{2})+ d_2
\sin(\lambda (2\pi-\theta))) &&\hbox{ for }
\frac{\pi}{2}<\theta<2{\pi},\eeqs where $\lambda\in (0,1)$ is given
by
\[
\lambda=\frac{2}{\pi}\arccos\left(\frac{1-\mu}{2|1+\mu|}\right),
\]
and the constants $c_1, c_2, d_1, d_2$ are appropriately defined.

\noindent Now we show using the arguments of Section \ref{sbvp} that
for $-\frac13<\mu<0$ and $\mu<-3$, the assumption
(\ref{tcoercivite}) holds. As before we define
\[
{\cal R}(v_+)(x,y)= \left\{\begin{array}{lll} v_+(-x,y)\quad
&\forall (x,y)\in (-1,0)\times (0,1),\\
v_+(-x,-y)\quad
&\forall (x,y)\in (-1,0)\times (-1,0),\\
v_+(x,-y)\quad &\forall (x,y)\in (0,1)\times (-1,0). \end{array}
\right.
\]
This extension defines an element of $H^1_-(\Omega_-)$ such that
\[
{\cal R}(v_+)= v_+ \quad \hbox{ on } \Sigma.
\]
Moreover with this choice we have
$$
\sup_{v\in H^1_+(\Omega_+)\atop v\neq 0}\frac{|{\it B}_-({\cal R}
(v),{\cal R} (v))|}{{\it B}_+(v,v)}=3|\mu|,
$$
and therefore for
\[
3|\mu|<1,
\]
we deduce   that (\ref{tcoercivite}) holds.
\\
To exchange the role of $\Omega_+$ and $\Omega_-$ we  define the
following extension from $\Omega_-$ to $\Omega_+$: for $v_-\in
H^1_-(\Omega_-)$, let
\[
{\cal R}(v_-)(x,y)=  v_-(-x,y)+v_-(x,-y)-v_-(-x,-y) \quad \forall
(x,y)\in \Omega_+.
\]
We readily check that it defines an element of $H^1_+(\Omega_+)$
such that
\[
{\cal R}(v_-)= v_- \quad \hbox{ on } \Sigma.
\]
Moreover with this choice we have (using the estimate $(a+b+c)^2\leq
3 (a^2+b^2+c^2)$ valid for all real numbers $a, b, c$)
\[
\sup_{v\in H^1_-(\Omega_-), v\ne 0}\frac{{\it B}_+({\cal R}
(v),{\cal R} (v))}{|{\it B}_-(v,v)|}\leq 3/|\mu|,
\]
and therefore for
\[
3/|\mu|<1,
\]
we deduce   that (\ref{tcoercivite}) holds.

For this second test, we take as  exact solution   the singular
function $u(x,y)=S(x,y)$  for $\mu=-5$ and  $\mu=-100$,
non-homogeneous Dirichlet boundary conditions on $\Gamma$ are fixed
accordingly. First, with uniform meshes, we obtain the expected
convergence rate of order $\lambda$ (resp. $2\lambda$) for the $
H^1$ (resp. $L^2$) error norm, see Tables \ref{table_sing_5_unif}
and \ref{table_sing_100_unif}. There, for sufficiently fine meshes,
we may notice that the effectivity index varies between 1 and 0.6
for $\mu=-5$ or between 9 and 6 for $\mu=-100$. From these results
we can say that the effectivity index depends on $\mu$, this is
confirmed by the numerical results obtained by an adaptive algorithm
(see below).
\begin{table}[ht]
\centering
\begin{tabular}{|l|l|c|c|c|c|c|}
\hline
&&&&&& \\
$k$ & $DoF$ & $e_{L^2}(u_h)$ & $CV_{L^2}$ & $ e_{H^1}(u_h)$ & $CV_{H^1}$ & $\displaystyle \frac{\eta(u_h)  }{e_{H^1}(u_h)  }$  \\
&&& &&&\\
\hline
1 & 289  & 1.60E-02 & & 2.84E-01& &2.57 \\
\hline
2 & 1089 & 8.66E-03 &0.93& 2.10E-01&0.45& 1.94 \\
\hline
3 & 4225 & 4.63E-03 &0.92& 1.55E-01& 0.45& 1.46\\
\hline
4 & 16641 & 2.47E-03 &0.92& 1.13E-01& 0.45 &1.09\\
\hline
5 & 32761 &1.80E-03 &0.92& 9.69E-02&0.46&0.95\\
\hline
6 & 90601 &1.13E-03 &0.92& 7.68E-02& 0.46&0.76\\
\hline
7 &  251001& 7.08E-04  & 0.92 & 6.08E-02  & 0.46 & 0.61\\
\hline
\end{tabular}
\caption{The singular solution, $\mu=-5 $, $\lambda \approx 0.46 $ (uniform refinement).}
\label{table_sing_5_unif}
\end{table}

\begin{table}[ht]
\centering
\begin{tabular}{|l|l|c|c|c|c|c|}
\hline
&&&&&& \\
$k$ & $DoF$ & $e_{L^2}(u_h)$ & $CV_{L^2}$ & $ e_{H^1}(u_h)$ & $CV_{H^1}$ & $\displaystyle \frac{\eta(u_h)  }{e_{H^1}(u_h)  }$  \\
&&& &&&\\
\hline
1 & 289  & 6.12E03 & & 1.54E-01& &18.77 \\
\hline
2 & 1089 & 2.59E-03 &1.29& 9.91E-02&0.66& 15.04 \\
\hline
3 & 4225 & 1.08E-03 &1.29& 6.35E-02& 0.66& 12.06  \\
\hline
4 & 16641 & 4.46E-04 &1.29& 4.04E-02&0.66& 9.66 \\
\hline
5 & 32761 &2.88E-04 &1.29& 3.24E-02&0.66& 8.65\\
\hline
6 & 90601 &1.49E-04 &1.30& 2.32E-02& 0.66 & 7.33\\
\hline
7 &  251001& 7.66E-05  & 1.30 & 1.66E-02 &0.66 & 6.21  \\
\hline
\end{tabular}
\caption{The singular solution, $\mu=-100 $, $\lambda \approx 0.66 $ (uniform refinement).}
\label{table_sing_100_unif}
\end{table}

Secondly, an adaptive mesh refinement strategy
 is used based on the estimator $\eta_T=\eta_{R,T}+\eta_{J,T}$, the marking procedure
 $$
 \eta_T>0.5 \max_{T'}\eta_{T'}$$
 and a standard refinement procedure with a limitation on the
 minimal angle.

For   $\mu=-5$ (resp. $\mu=-100$),  Table
\ref{table_dis_sing_a5_adap} (resp. \ref{table_dis_sing_a100_adap})
displays the  same quantitative results as before. There we see that
the
 effectivity index is around 3 (resp. 34), which is quite satisfactory
 and comparable with results from \cite{cochez:10,ern:07}.
 As before and in these references we notice that it deteriorates as the contrast
 becomes
 larger. On these tables we also remark a convergence order of 0.76
 (resp. 1) in the $H^1$-norm and mainly the double in the
 $L^2$-norm. This yields   better orders of convergence as for
 uniform meshes as expected, the case $\mu=-5$ giving less accurate
 results due to the high singular behavior of the solution (a similar phenomenon occurs in \cite{cochez:10}  for
 instance).

\begin{table}[ht]
\centering
\begin{tabular}{|l|l|c|c|c|c|c|c|c|}
\hline
&&&&&&\\
$k$ & $DoF$ & $e_{L^2}(u_h)$ & $CV_{L^2}$ & $ e_{H^1}(u_h)$ & $CV_{H^1}$ & $\displaystyle \frac{\eta(u_h)  }{e_{H^1}(u_h)  }$  \\
&&&&&&\\
\hline
1 & 81 & 2.92E-02& & 3.79E-01 & &3.39 \\
\hline
5 & 432 & 3.49E-03 &2.54& 1.40E-01 &1.19 & 4.18 \\
\hline
7 & 1672 & 1.25E-03 &1.52& 8.04E-02 &0.82 & 4.07 \\
\hline
10 & 5136 & 4.26E-04 &1.92& 4.90E-02 &0.88 & 3.63 \\
\hline
13 & 20588 & 1.64E-04  &1.37&3.14E-02 &0.64& 3.32  \\
\hline
18&80793&5.50E-05&1.60&1.80E-02&0.81&3.23\\
\hline
24&272923&2.39E-05&1.37&1.17E-02&0.71&2.5\\
\hline
\end{tabular}
\caption{The singular solution, $\mu=-5 $, $\lambda \approx 0.46 $
(local refinement).} \label{table_dis_sing_a5_adap}
\end{table}

\begin{table}[ht]
\centering
\begin{tabular}{|l|l|c|c|c|c|c|c|}
\hline
&&&&&& \\
$k$ & $DoF$ & $e_{L^2}(u_h)$ & $CV_{L^2}$ & $ e_{H^1}(u_h)$ & $CV_{H^1}$ & $\displaystyle \frac{\eta(u_h)  }{e_{H^1}(u_h)  }$  \\
&&&&&&\\
\hline
1 & 81 & 1.41E-02& & 2.35E-01 & & 23.59  \\
\hline
4 & 363 & 1.93E-03 &2.65& 8.77E-02 &1.31 & 34.86 \\
\hline
7 & 1566 & 4.94E-04 &1.86& 4.31E-02 &0.97& 33.10 \\
\hline
11 & 5981 & 1.23E-04  &2.07&2.15E-02 &1.04& 33.17  \\
\hline
16&25452&2.98E-05&1.96&1.05E-02&0.99&34.65\\
\hline
24&106827&7.36E-06&1.95&5.23E-03&0.97&33.89\\
\hline
\end{tabular}
\caption{The singular solution, $\mu=-100 $, $\lambda \approx 0.66 $
(local refinement). \label{table_dis_sing_a100_adap}}
\end{table}

\clearpage

\newcommand{\noopsort}[1]{}\def\cprime{$'$}

\end{document}